\documentclass{article}
\usepackage{hyperref}

\def\E{\mathop{\hbox{\sf E}}\nolimits}
\def\P{\mathop{\hbox{\sf P}}\nolimits}

\def\dd{\displaystyle}
\def\nn{\nonumber}

\def\phi{\varphi}

\def\f{\frac}
\def\d{\partial}
\def\G{\Gamma}

\def\B{\Big}

\def\bs{\backslash}
\def\bB{{\bf B}}
\def\La{\Lambda}

\newtheorem{lemma}{Lemma}[section]
\newtheorem{theorem}{Theorem}
\newtheorem{definition}{Definition}[section]

\newtheorem{conjecture}{Conjecture}

\usepackage{graphics}

\def\figref#1{fig. \ref{fig.#1}}
\def\eqref#1{(\ref{#1})}

\def\putfigure#1#2{
	\begin{figure}[ht]
	\centering
	\includegraphics{#1.eps}
	\caption{#2}
	\label{fig.#1}
	\end{figure}
}

\def\defined#1{{\em #1}}

\def\optional#1{}

\def\EQ#1#2{\begin{equation}\label{#1} #2 \end{equation}}

\def\skel{\mbox{skel}}
\def\mod{\mbox{\rm\,mod\,}}

\def\bibauthor{}
\def\bibtitle{}
\def\bibjournal{}

\def\prodl{\prod\limits}
\def\suml{\sum\limits}
\def\cupl{\bigcup\limits}

\def\proof{\par\noindent{\bf Proof.\ }}
\def\proofof#1{\par\noindent{\bf Proof of #1.\ }}
\def\eop{\vskip 3mm }

\def\RN{\stackrel{rn}{\subset}}
\def\XX{F_{anc}^{(r)}}

\def\bM{{\bf M}}
\def\bA{{\bf A}}

\def\calL{{\cal L}}
\def\calT{{\cal T}}

\def\barB{{\bar B}}
\def\barphi{{\bar \phi}}
\def\baromega{{\bar \omega}}
\def\barzeta{{\bar \zeta}}

\begin{document}
%\title{Exact results for the UIPT}
\title{Uniform infinite planar triangulation and related time-reversed critical branching process%
    \footnote{This research was partially supported by RFBR grant 02-01-00415.}}%
\author{Maxim Krikun%
    \thanks{
    Laboratory of Large Random Systems,
    Faculty of Mechanics and Mathematics,
    Moscow State University.
    E-mail: krikun@lbss.math.msu.su}}

\maketitle

\begin{abstract}
We establish a connection between the uniform infinite planar triangulation
and some critical time-reversed branching process. This allows to find 
a scaling limit for the principal boundary component of a ball of radius $R$ for large $R$ 
(i.e. for a boundary component separating the ball from infinity).
We show also that outside of $R$-ball a contour exists that has length linear in $R$.
\end{abstract}

\section*{Introduction}

%This work is based on two recent papers by Angel and Schramm.
The \defined{uniform infinite planar triangulation (UIPT)} is a random graph, considered 
as one of possible models of {\em generic planar geometry}.
UIPT is defined as a weak limit of uniform measures on triangulations with finite number of triangles.
%\footnote{Similar models were considered before in the literature, however this is \cite{AS} where
%the existence of UIPT is proved.}
In \cite{AS} Angel ans Schramm proved the existence of this limit,
in \cite{A} some basic geometrical properties of UIPT were investigated, in particular it was found 
that the ball of radius $R$ in UIPT has volume of order $R^4$ and boundary of order $R^2$, 
up to polylogarithmic terms.
This fact reflects the conjecture known in physics, see \cite{Amb}.

In this paper we improve the result of \cite{A} concerning the boundary of a ball and
give an exact limit of a corresponding scaled random variable as $R\to\infty$.
We use a new combinatorial "skeleton" construction, which uncovers a connection between UIPT profile 
and certain time-reversed branching process. Using this connection we state a new fact
concerning the UIPT: we show that outside of the $R$-ball a contour exists, that 
separates the ball from the infinite part of triangulation and has length linear in $R$.

The paper is organized as follows. In the first part we give necessary definitions and 
review some results of \cite{AS,A}.  
The main results of this paper are stated in section \ref{section.results}. 
In the second part we describe the skeleton construction.
In the third part we use the "raw method" to obtain the limiting distribution of $R$-ball 
boundary length as $R\to\infty$.
In the fourth part we show how the UIPT is related to time-reversed branching processes and
derive the existence of linear contour.
In the last part we discuss the universality of the obtained results abd show briefly
how the techniques developed in this work can be applied to some different flavor of 
triangulations. Some technical proofs are contained in the appendix.  
No historical review of the subject is provided in this paper and only the results 
which are necessary for the statement of the problem and proof of the main theorems 
are mentioned.

\section{Definitions and results}

\subsection{Rooted near-triangulations}

We'll start from the very beginning.
Planar maps and planar triangulations were studied from the 60-ths, starting with works of 
Tutte \cite{Tutte}. The following definitions (\ref{def.map}--\ref{def.rnt}) are given here
according to \cite{GJ}.
\begin{definition}\label{def.map}
A \defined{planar map} $M$ is a nonempty connected graph $G$ embedded into a sphere $S$.
It separates the sphere into disjoint areas, called \defined{faces}. The number of edges
incident to a face is called the \defined{face degree}. 
%Two maps are equivalent if a homeomorphism of a sphere exists that maps one to another.
\end{definition}
When dealing with the combinatorics of maps the usual practice is to consider rooted maps 
instead of generic ones in order to avoid problems with non-trivial symmetries.
It's known that almost all maps in some classes have no nontrivial authomorphisms (\cite{RW}),
so adding a root hopefully will not affect the results in most cases.  
\begin{definition}\label{def.rmap}
A \defined{root} in a map consists of a face and a directed edge on 
the boundary of that face, which are called \defined{rooted face} and \defined{root edge} 
respectively. The tail vertex of the rooted edge is called a \defined{root vertex}.
A \defined{rooted planar map} is a planar map with a root.
%The sphere homeomorphism mapping to equivalent rooted maps one to another should respect the root.
\end{definition}

\begin{definition}\label{def.cutvertex}
A vertex $v$ in a connected graph $G$ is called a \defined{cut vertex}, if $G$ can be 
decomposed in two connected graphs $G_1$ and $G_2$, both containing at least one edge,
such that their intersection consists of a single vertex $v$.
%$G_1 \cup G_2 = G$, $G_1 \cap G_2 = \{v\}$.
\end{definition}

\begin{definition}\label{def.rnt}
A \defined{rooted near-triangulation (RNT)} is a rooted planar map, such that all it's faces 
except the rooted face are triangles (i.e. have degree $3$) and it contains no cut vertices.
\end{definition}
In \cite{AS} the rooted near-triangulations are called \defined{type II triangulations},
and by \defined{type III triangulations} is denoted a class of rooted near-triangulations 
with no multiple edges.
We'll consider type III triangulations in the last part, when discussing the 
universality conjecture.

The main result concerning RNT we'll need is the following:
the number of RNTs with $N$ triangles and $m$ edges on the boundary is
\begin{equation}\label{eq.tutte}
C_0(N,m) = \f{2^{j+2}(2m+3j-1)!(2m-3)!}{(j+1)(2m+2j)!((m-2)!)^2}, \qquad N=m+2j
\end{equation}

\subsection{UIPT}

In \cite{AS} the infinite triangulations are defined identically to the finite ones, 
except that they are required to be \defined{locally finite}.
\begin{definition}\label{def.locfin}
A triangulation $T$ embedded into a sphere $S$ is \defined{locally finite}, if 
every point in $T$ has a neighborhood that intersects only a finite number of elements
of $T$ (i.e. edges, vertices, triangles).
\end{definition}
In particular this means that each vertex of $T$ has finite degree.

Let $\calT$ be the space of finite and infinite triangulations with a natural metrical 
topology: the distance between two triangulations $T_1$ and $T_2$ is $(k+1)^{-1}$, where 
$k$ is the maximal radius such that two balls around the root in $T_1$ and $T_2$ are
equivalent. This topology induces weak topology in the space of measures on $\calT$,
and weak convergence of (probability) measures on $\tau$. Namely
\begin{definition}
A measure $\tau$ on $\calT$ is a limit of $\tau_n$ if for every bounded continuous function
$f: \calT \to R$ 
\[ \lim_{n\to\infty} \int f d \tau_n = \int f d \tau. \]
\end{definition}
Now let $\tau_n$ be a uniform probability measure on the set of triangulations with $n$ triangles.
\begin{theorem}[Angel, Schramm]\label{thm.as}
There exists the probability measure $\tau$ supported on infinite planar triangulations,
such that 
\[ \lim \tau_n = \tau. \]
\end{theorem}
Moreover, for theorem \ref{thm.as} to hold it's enough to show that a probability measure $\tau$
exists, such that for every radius $R$ and every triangulation $T$
\begin{equation}\label{eq.limtau}
\tau_n(B_R=T) \to \tau(B_R=T).
\end{equation}

In fact, Theorem \ref{thm.as} can be stated not for the class of full sphere triangulations, but 
for a class of rooted near-triangulations with fixed boundary length, i.e. $\tau_n(m_0)$ being 
a uniform distribution on RNTs with $n$ triangles and boundary $m_0$, there should exists a 
limit $\tau(m_0)$, which is a uniform distribution of infinite triangulations with rooted 
boundary of length $m_0$. The existence of such a limit follows from the observation that the
distribution $\tau_n(m_0)$ can be considered as $\tau_{n+m_0-1}$ conditioned to have 
$m_0$ triangles around the root vertex. 

By $S(N)$ and $S_\infty$ we shall denote samples of measures $\tau_N$ and $\tau$ respectively,
and by $S(N,m_0)$, $S_\infty(m_0)$ the samples of $\tau_N(m_0)$, $\tau(m_0)$.

\subsection{Multi-rooted triangulations}

One can consider RNT as a triangulation of a disk, or of a sphere with a hole. 
The disk is obtained from a sphere by cutting the rooted face of a RNT.
We shall generalize the definition of RNT to include triangulations of a sphere with multiple 
holes.
\begin{definition}\label{def.mrt}
A \defined{multirooted triangulation (MRT)} of type $(N,m_0,k;m_1,\ldots,m_k)$ is 
a rooted planar map, such that
\begin{itemize}
\item  the rooted face has has degree $m_0$,
\item  $k$ faces are distinguished and labeled with numbers $1,\ldots,k$, these faces 
       are called holes.
\item  the holes have degrees $m_1,\ldots,m_k$, 
\item  on the boundary of each hole a directed edge (additional root) is specified,
       so that it's orientation coincides with the orientation of the rooted face,
\item  there are $N$ more triangular faces,
\item  there is no cut vertices. 
\end{itemize}
%The sphere homeomorphism mapping two equivalent multirooted triangulations one to another 
%should respect holes labeling and additional roots.
\end{definition}
In the following two parts of the paper we will refer to MRT simply as triangulation.
We will also use the notation $N(T)$, $m_0(T)$, $m_j(T)$ to denote the parameters 
of a triangulation.
In order to easily refer a particular vertex on the boundary, we use the 
\defined{\label{def.stdenum}standard boundary enumeration}: 
we enumerate the vertices on the boundary
in clockwise direction, starting from number $0$ for the root vertex, so that the 
root edge starts at vertex $0$ and ends at vertex $1$.

We do not impose any restrictive conditions when defining MRT: 
two boundaries of a MRT can share and edge, there can even be no internal triangular 
faces at all ($N(T)=0$). 
The following two definitions outline a useful subclass of rigid triangulations.

\begin{definition}\label{def.completion}
Given a triangulation $T$ of type $(N,m_0;m_1,\ldots,m_k)$, $k\ge1$, call a sequence of 
disk triangulations $D_1,\ldots,D_k$ an \defined{appropriate disk set}, if 
$m_0(D_j)=m_j$, $j=1,\ldots,k$.

Glue each disk $D_j$ to a corresponding boundary $m_j$ of $T$, so that the root of $D_j$ 
coincides with a rooted edge on $m_j$. The result is a triangulation $T'$ of type
$(N',m_0)$, $N'=N+\sum_j N(D_j)$. Call this operation a \defined{completion to 
the sphere} and denote it by an equation $T'=T+(D_1,\ldots,D_k)$.
\end{definition}

\begin{definition}\label{def.rigid}
A triangulation $T$ is \defined{rigid}, if for distinct appropriate disk sets 
$D_1,\ldots,D_k\neq D'_1,\ldots,D'_k$,
the result of completion differs: $T+(D_1,\ldots,D_k) \neq T+(D'_1,\ldots,D'_k)$,
\end{definition}

Rigidity is an essential property for counting triangulations and subtriangulations.
If $T$ is rigid, it can be a subtriangulation of $S$ in one only position.
\begin{definition}\label{def.rn}
A rigid triangulation $T$ of type $(N,m_0,k;m_1,\ldots,m_k)$ is a \defined{root neighborhood} 
in a triangulation $S$ of type $(N,m_0,0;)$, if $T$ can be completed to $S$ with an appropriate set of 
disks. %, i.e. $S(N,m)=T+(D_1,\ldots,D_k)$.
Denote this by $T\RN S$.
\end{definition}

Let us return to the formula \eqref{eq.tutte}. It describes the number of triangulations
of type $(N,m,0)$. The corresponding generating function is 
\begin{equation}\label{eq.u0}
U_0(x,y) = \sum_{N=0}^\infty \sum_{m=2}^\infty C_0(N,m) x^N y^{m-2} 
  = \f{ y-x+(h-y)\sqrt{\f{x^2}{h^2}-4xy} }{2xy^2},
\end{equation}
where the function $h=h(x)$ is the solution of \EQ{eq.defh}{ x = h-2h^3 } such that $h(0)=0$
.

Most probabilities in this paper arise from singularity analysis of this function 
(for details of $U_0(x,y)$ analysis see also \cite{MK}).
Fix $y$ such that $|y|<y_0=\f1{\sqrt6}$, then the principal singularities of $U_0(x,y)$
as a function of $x$ are the two points $x=\pm x_0$, $x_0=\sqrt{2/27}$. Near these
points an expansion holds
\[ U_0(x,y)\B|_{x= x_0-t} = A(y)+A_1(y)t+\f{6^{3/4}}3 B(y)t^{3/2}+O(t^2), \]
\[ U_0(x,y)\B|_{x=-x_0+t} = A(-y)+A_1(-y)t+\f{6^{3/4}}3 B(-y)t^{3/2}+O(t^2), \]
where 
\[ A(y) = \f32\f{2\sqrt{1-\sqrt6y}+1}{(\sqrt{1-\sqrt6y}+1)^2}, \qquad
   B(y) = \f{1}{(1-\sqrt6y)^{3/2}}, \]
and $A_1(y)$ does not play any role in further calculations. Let
\[ a(m) = [y^{m-2}] A(y), \qquad
   b(m) = [y^{m-2}] B(y). \]
We now state two theorems concerning the limiting distribution $\tau$. 
These appear in \cite{AS} in slightly different notation (and are referred there as 
not entirely new). The proof using $U_0(x,y)$ singularity analysis is rather 
straightforward, so we leave it to the appendix. 

%Let $S(N,m)$ denote a sample of uniform measure on triangulations of type $(N,m,0)$.
\begin{theorem}\label{thm.limdist}
Given a rigid triangulation $T$ of type $(n,m_0,k;m_1,\ldots,m_k)$, there exists a limit
\begin{equation}\label{eq.limdist}
\lim_{N\to\infty} \P\{ T\RN S(N,m_0) \} 
   = \f1{b(m_0)} a(m_1)\cdots a(m_k) \sum_{j=1}^k \f{b(m_j)}{a(m_j)} x_0^n. 
\end{equation}
\end{theorem}
Further we will denote the limit \eqref{eq.limdist} by $\P\{T\RN S_\infty(m_0)\}$.
Note that in a particular case $k=1$, when $T$ has type $(n,m_0,1;m_1)$,
\eqref{eq.limdist} becomes
\begin{equation}\label{eq.limdist2}
\lim_{N\to\infty} \P\{ T\RN S(N,m_0) \} = \f{b(m_1)}{b(m_0)} x_0^n. 
\end{equation}

\begin{theorem}\label{thm.oneend}
Given a triangulation $T$ of type $(n,m_0,k;m_1,\ldots,m_k)$, 
consider random triangulation $S(N,m_0)$ 
under condition $T\RN S(N,m_0)$. Let $S(N,m_0)=T+(D_1,\ldots,D_k)$ and let $N_j=N(D_j)$.
Take a limit $N\to\infty$. Then the largest of $N$s is infinite while others are a.s. finite:
\[ \lim_{N\to\infty} \P\{ N_j=\max_{i=1,\ldots,k}(N_i) \} = \f{b(m_j)}{a(m_j)}\B(\sum_{i=1}^k
\f{b(m_i)}{a(m_i)}\B)^{-1}, \]
and the limiting conditional distribution exists
\[ \lim_{N\to\infty}\P\{ N_i=n_i, i=1,\ldots,\widehat j,\ldots,k | N_j=\max(N_i) \}. \]
\end{theorem}

\subsection{Main results}\label{section.results}
Main results of this paper are summarized in this section.
(Note, that we use here the notion of $R$-hull, which is defined in section \ref{section.ball} below.
The $R$-hull is a natural modification of $R$-ball, such that it's boundary always consists of 
a single component. One can think of it as a {\em disk of radius $R$ centered in the root}).
\begin{theorem}\label{thm.profile}
The upper boundary $m_1(\barB_R)$ of a $R$-hull of UIPT (with $m_0=2$) grows as $R^2$. 
There exists a limit
\[ \lim_{R\to\infty} \f{m_1(\barB_R)}{R^2} = \xi, \]
where $\xi$ is a random variable with density
\[ p_{\xi}(x) = \f{2}{\sqrt\pi}e^{-t}t^{1/2}. \]
\end{theorem}

\begin{theorem}\label{thm.bp}
Let $\barB_R$ be an $R$-hull of UIPT. 
\[
\P\{ m_1(B_R)=k | m_0(B_R)=l \} 
  %= \f{b(k)y_0^{-k}/k}{b(l)y_0^{-l}/l} \P\{ \barzeta(R)=l | \barzeta(0)=k \}, 
   = \f{[t^k] F_0(t)}{[t^l] F_0(t)} \P\{ \barzeta(R)=l | \barzeta(0)=k \}, 
\]
where
\[ F_0(t) = \f2{\sqrt{1-t}} -4 + 2\sqrt{1-t}. \]
and $\barzeta$ is a critical branching process with special behavior near zero (see section 
\ref{section.mbp} for details).
\end{theorem}

\begin{theorem}\label{thm.lc}
For each $R$ there exists a (random) contour in UIPT such that 
\begin{itemize}
\item it lies outside of $\barB_R$,
\item it separates root from infinity,
\item it's expected length is linear in $R$ as $R\to\infty$.
\end{itemize}
\end{theorem}

\section{UIPT representation}

\subsection{The ball}\label{section.ball}
We use a simple combinatorial metric on triangulations, where each edge has length one and
the distance between two points equals the number of edges in the shortest path between them.
To each vertex we assign it's height $H(v)$, which is a minimal distance to the rooted boundary.
Obviously, the value of $H$ on the ends of an edge can differ at most at one, so
each triangle matches one of the three patterns: $(R,R,R+1)$, $(R,R,R-1)$, $(R,R,R)$. 
We call such triangles \defined{plus-}, \defined{minus-} and \defined{zero-triangles}
respectively.

\begin{definition}\label{def.ball}
The \defined{ball} of radius $R$ $B_R$ consists of all triangles (including edges and vertices of
a triangle) that have at least one vertex with $H(v)\le R-1$.
\end{definition}
Being defined in such a way, $B_R$ contains all vertices with height $H(v)=R$ but not all 
edges with both ends at height $R$. 
For such an edge to be included into $B_R$ it should be an upper edge of a minus-triangle
with vertices heights $(R,R,R-1)$. The reason of such definition choice is the following:
\begin{lemma}
Take a sphere triangulation $S$ and a ball $B_R=B_R(S)$. 
Take some completion $S'=B_R+(D_1,\ldots,D_k)$, so that $S'\neq S$. Then the ball remains
the same, $B_R(S')=B_R(S)$.

In other words, two distinct balls of the same radius, $B_R$ and $B_R'$, at most one of them
can be a root neighborhood in a triangulation $S$.
\end{lemma}
This would not be true, if all the $(R,R)$ edges would be included into the ball. 

The ball boundary is not necessary connected, and $S\bs B_R$ may consist of several disjoint parts,
$S\bs B_R = D_1, \ldots, D_k$. However, due to the Theorem \ref{thm.oneend}, in UIPT 
a.s. only one of disks $D_1,\ldots,D_k$ contains infinitely many triangles. 
If we'll cut this disk only and keep all others, we'll get an a.s. finite root neighborhood 
with a single boundary, where all the vertices of the boundary have the same height $R$.
Such root neighborhood is called a \defined{hull}.
\begin{definition}\label{def.hull}
Consider a triangulation $T$ and a ball $B_R$. The set $T\bs B_R$ may not be connected. Generally
this set consists of disk triangulations $D_1,\ldots,D_k$. Let $D_j$ be a disk with maximal number 
of triangles. Then the \defined{hull} $\barB_R$ consists of all the triangles contained in $T$, 
but not in $D_j$: $\barB_R = T\bs D_j$.
\end{definition}
Following the definition \ref{def.mrt}, $\barB_R$ is a triangulation of type $(N,m_0,1;m_1)$.
Let us call such triangulations \defined{cylindric}. For cylindric triangulations we will
refer to $m_0$ and $m_1$ as \defined{lower} and \defined{upper} boundaries respectively.

\subsection{Skeleton construction}
Now consider a UIPT sample $S$ and a sequence of increasing hulls $B_1,\ldots,B_R$ in $S$.
As a consequence of Theorem \ref{thm.oneend}, $\bar\barB_1\subset \ldots \subset \barB_R$
(note that this would be not true in general for finite $S$).

\begin{definition}\label{def.layer}
A \defined{layer} $L_R$ is a subtriangulation in $S$ that consists of triangles contained
in $\barB_R$ but not in $\barB_{R-1}$.\footnote{
  Note that this definition is incomplete, since definition \ref{def.mrt} requires 
  an additional root on the upper boundary to be defined.  }
\end{definition}
Let $\calL_R$ be a mapping from the subset $\calT_\infty$ of infinite triangulations in$\calT$
to cylindric triangulations $\calL_R: S \to L_R(S)$, and let 
$\La = \cupl_{R=1}^\infty \calL_R(\calT_\infty)$ be the set of all possible layers.
Further we'll show that $\calL_R(\calT_\infty)=\La$ for all $R$, 
so each layer $L\in\La$ may appear at any level $R$ within some triangulation.

Let us examine some properties of a layer $L\in\La$:
\begin{itemize}
\item $L$ is a cylindric triangulation with all vertices of the upper boundary having height $1$;
\item each edge of the upper boundary of $L$ belongs to some minus-triangle, which has one vertex 
      at the lower boundary of $L$
\item between two such subsequent minus-triangles a disk triangulation is contained, which has
      one vertex at the upper boundary and at least one vertex at the lower boundary. 
      %(if two subsequent minus-triangles do share an edge, say that a triangulation $\delta$ 
      %of type $(0,2,0)$ is contained between them).  
\end{itemize}
These properties follow immediately from the definition of a layer.
\begin{definition}\label{def.lp}
Given a layer $L$, the corresponding \defined{layer pattern (LP)} is a subgraph in $L$ that 
consists of it's upper boundary, lower boundary and all minus-triangles with an edge on the 
upper boundary of $L$. 

The boundaries between two subsequent minus-triangles are called
\defined{slots} of a LP.
For convenience, if two minus-triangles in a LP share an edge, we split this
edge and put a slot of length two, so that each triangle in LP is followed by a slot.

A slot and a triangle to the right of it are called \defined{associated}. 
The additional root on the slot boundary is the common edge of a slot and an associated 
minus-triangle.
The additional root on the upper boundary of a layer is the upper edge of a minus-triangle,
associated to the slot that contains the root of the lower boundary.\footnote{
    the same definition should be used for the additional root on the upper boundary 
    of a layer to complete definition \ref{def.layer}.                       }
\end{definition}
\putfigure{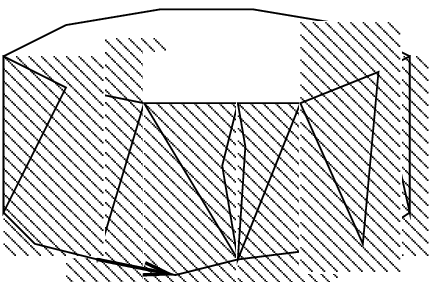}{A layer pattern with two 4-slots, one 2-slot and two 3-slots visible. }
Each layer matches exactly one layer pattern. A layer can be obtained from the corresponding LP
by filling the slots with an appropriate disk set (this operation is similar to definition
\ref{def.completion}, except that the upper boundary is kept open), and vice versa, 
almost every appropriate disk set can be used to obtain a layer from the layer pattern.
There is however one exception.
\begin{lemma} 
Given a layer pattern $P$ with an upper boundary of length $k$ and $k$ slots of length
$l_1,\ldots,l_k$ and an appropriate disk set $D_1,\ldots,D_k$, (i.e. such that the 
boundary length of a disk corresponds to the length of a slot, $m_0(D_j)=l_j$).

Glue each disk to the corresponding slot, so that the root of a disk matches the common
edge of a slot and it's associated triangle. Then the result is a layer $L\in\La$, 
except that 
$$
\begin{array}{ll}
\mbox{1)}&\mbox{all of $l_1,\ldots,l_k$ are equal to $2$, except $l_j$;          } \\
\mbox{2)}&\mbox{the disk $D_j$ contains an edge between the vertices $0$ and $2$ }\\
         &\mbox{of the rooted boundary.}
\end{array}
\eqno (E)
$$
(see boundary enumeration defined at page \pageref{def.stdenum}).
\end{lemma}
\putfigure{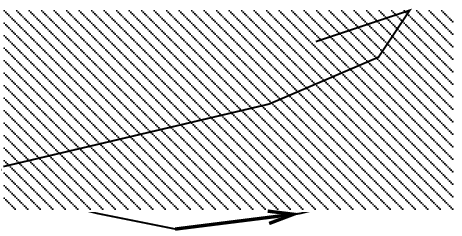}{A layer pattern that falls under exception (E).}
\proof
First let us check that $L=P+(D_1,\ldots,D_k)$ confirms MRT definition.
Only the last condition in Definition \ref{def.mrt} (the one that requires no cut 
vertices in triangulation) is nontrivial.

In a planar triangulation with multiple boundaries a cut vertex may appear in two ways:
either some boundary is \defined{self-touching}, i.e. some vertex is met twice when walking around 
this boundary, either there is a loop -- an edge with both ends at the same vertex. 

The upper and lower boundaries of LP are not self-touching. When gluing disks
$D_1,\ldots,D_k$ no vertices of LP are glued together, so the boundaries remain 
not self-touching. Hence we should check for loops only.

A loop may appear if some slot boundary is self-touching in some vertex $v=v'$, and this 
boundary is filled with a disk so that the ends of some edge in a disk are identified 
with $v$ and $v'$; this is exactly the case described by exception (E).
\eop

\begin{lemma}\label{lemma.lseq}
Given a sequence of layers $L_1,\ldots,L_R$, such that $k(L_i)=l(L_{i+1})$, by gluing 
upper boundary of $L_i$ with lower boundary of $L_{i+1}$ for $i=1,\ldots,R-1$ we get
a valid triangulation; consequently the image of $\calL_R$ is the same for all $R$:
\[ \calL_1(\calT_\infty)=\ldots=\calL_R(\calT_\infty)=\La. \]
\end{lemma}
\proof The proof is similar to that of the previous lemma. We have to check only that no cut vertices
appear when adding a layer $L_j$ to a cylindric triangulation $T_{j-1}=L_1+\ldots+L_{j-1}$.
This is true, since both $L_j$ and $T_{j-1}$ contain no cut vertices, and they are glued at least in 
two points, so the graph $T_j=T_{j-1}+L_j$ also has no cut vertices.
\eop

\begin{definition}\label{def.skel}
Given a hull $\barB_R$ in a triangulation $T$, it's \defined{skeleton} $\skel(B_R)$
is a subtriangulation in $\barB_R$
that contains all the minus-triangles with an upper edge on the boundary of $\barB_1(T),\ldots,\barB_R(T)$.

Let $L_1,\ldots,L_R$ is sequence of layers corresponding to $\barB_R$ and $P_1,\ldots,P_R$ 
a sequence of layer patterns ($L_j$ matches $P_j$ for $j=1,\ldots,R$), then
\[ \skel(\barB_R) = P_1 + \ldots + P_k. \]  
The skeleton is a multi-rooted triangulation, with additional roots specified by the layer patterns.
\end{definition}
There is a correspondence between skeletons and trajectories of branching processes (see \figref{skel}).
Given a skeleton $skel(\barB_R)$, say that
\begin{itemize}
\item a contour between two layers is a generation of a branching process;
\item the edges in this contour are particles;
\item the bottom edges of a slot are descendants of the upper edge in an associated triangle.
\end{itemize}
In section \ref{section.bp} we will explore this parallelism in details.
\putfigure{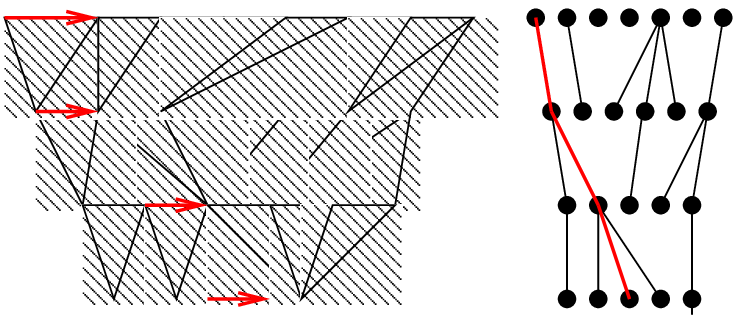}{A skeleton and a trajectory of a branching process}

\subsection{Skeleton enumeration}

The skeleton construction and Theorem \ref{thm.limdist} allow us to compute the probability
that an $R$-hull $\barB_R$ in an infinite random triangulation $S$ has a particular skeleton $K$. 
To do so, we shall take the set $\beta(K) = \{ T : \skel(B)=K \}$, assign to each triangulation it's
weight according to \eqref{eq.limdist2} and sum over $\beta(K)$,
\[ Z(K) = \sum_{T: \skel(T)=K} \f{b(m_1)}{b(m_0)} x_0^{n(T)}. \]
In fact, we could do this in a more general setting by replacing $x_0$ with $x$.

Since $m_0$ and $m_1$ are determined by the skeleton $K$, and hence are the same for all $T\in\skel^{-1}(K)$,
\[ Z(K) = \f{b(m_1)}{b(m_0)} \sum_{T:\skel(T)=K} x_0^{n(T)}. \]
Next, each $T\in\beta(K)$ is obtained from the skeleton $K$ by filling the slots $K$ with some disks,
and the disks in different slots are chosen independently.

Let the slots of $K$ have lengths $l_1,\ldots,l_q$. If we don't take the exception (E) into
account, (a.e. if $K$ doesn't contain layer patterns that fall under (E)),
each slot can be filled with any disk with appropriate boundary length. The sum $Z(K)$ then can 
be represented as a product over slots, where a slot $l_j$ has a term
\[ \sum_{n=0}^\infty C_0(n,l_j) x_0^{n+1} = x_0 [y^{l_j-2}] U_0(x_0,y), \]
an additional $x_0$ corresponds to the associated triangle. So 
\begin{equation}\label{eq.Zprod}
Z(K) = \f{b(m_1)}{b(m_0)} \prod_{j=1}^q x_0 [y^{l_j-2}] U_0(x_0,y). 
\end{equation}
In a more general case, when (E) is taken into account, some terms in the sum above 
should be modified, which is not a simple task for a generic skeleton.
In the following two section we do this in two different ways.

In section \ref{section.raw} we construct a sum, similar to $Z(K)$, that allows to 
enumerate all layers, then extend this result to generic $R$-hulls and finally obtain 
exact asymptotic for the upper boundary.

In section \ref{section.bp} we consider in details the relation between skeletons and 
branching processes. This way we get much simpler computations but no exact limits.
The main result in this section is the existence of a linear contour. 

\section{Raw approach}\label{section.raw}

\subsection{Layer statistical sum as a linear operator}
The sum $Z(K)$ does not take exception (E) into account.
However we can use the same argument to compute a statistical sum that enumerates all 
layers in $\La$ with respect to (the length of) both boundaries and the number of triangles.
\begin{lemma}\label{lemma.a}
Let $l=l(L)$ and $k=k(L)$ denote the length of lower and upper boundaries of a layer $L$,
and let $n=n(L)$ denote the number of triangles in $L$. 
Take a function $F(t)$ that allows expansion \[ F(t)=\suml_{j=0}^\infty c_j t^j. \]
Then we can write down the sum over all layers $L\in\La$
\[ \sum_L c_{k-2} x^n y^{l-2} = (\bA(x)F)(y), \]
where $\bA(x)$ is a linear operator
\begin{eqnarray}
(\bA(x)F)(y) &=& \f1y u(x,y) u'(x,y) F(u(x,y)) \nn\\
    &&{} - \f1y \B(u(x,y)+y u'(x,y)\B) u'(x,0) F(u(x,0)).
\label{Adef}
\end{eqnarray}
and $u(x,y) = x U_0(x,y)$,  $u'(x,y) = \f{\d}{\d y} x U_0(x,y)$.
\end{lemma}
\proof The proof is similar to the considerations used above for $Z(K)$, but there are two 
important things to note. 

Let $P$ be a layer pattern with upper boundary $k$ and lower boundary $l$, 
and let us enumerate it's slots, starting from the one that contains the root of lower boundary.
Let the slots have lengths $l_1,\ldots,l_k$. Then $(l_1-2)+\ldots+(l_k-2)=l$, and  
according to the definition of layer pattern, $l_1\ge3$, since the first slot should have 
at least one edge on the lower boundary. 
Given $l_1,\ldots,l_k$, the root can be placed in any of $(l_1-2)$ positions, so in order to 
define a layer pattern such position should be specified along with $l_1,\ldots,l_k$.

When translated to the generating functions language, this gives the first term in \eqref{Adef}:
\begin{eqnarray}
\sum_L c_{k-2} x^n y^{l-2} 
&=& y^{-2} \sum_{k=2}^\infty \sum_{l_1\ge3; l_2,\ldots,l_k\ge2} 
      c_{k-2} (l_1-2) y^{l_1-2} [t^{l_1-2}] u(x,t)        \nn\\
&&    %\phantom{\sum_{k=2}^\infty \sum_{l_1,\ldots,l_k\ge2} }
      \times \prod_{j=2}^k y^{l_j-2} [t^{l_j-2}] u(x,t)   \nn\\
&=& y^{-1} u'(x,y) u(x,y) F(u(x,y))
\label{Adef1}
\end{eqnarray}

In fact \eqref{Adef1} is not correct, since the term on the right includes some layers that should not be 
included due to exception (E), so we have to subtract something from \eqref{Adef1} to obtain a 
valid expression.

For a layer pattern $P$ to fall under exception (E) it should have $l_2=\ldots=l_k=2$,
for a corresponding layer $L$ to fall under (E) the disk triangulation glued into the first 
slot ($l_1$) should have an edge between the vertices $0$ and $2$ of the boundary (a $02$-edge). 
To count such triangulations we use the following statement.
\begin{lemma}\label{lemma.d}
Let $D(N,m,l)$ be the number of disk triangulations of type $(N,m+l,0)$ that have an edge between
the vertices $0$ and $l$ of the boundary (this may be a boundary edge too, when $m=1$ or $l=1$). 
Then
\begin{equation}\label{eq.d}
W(x,y,z) = \sum_{N\ge0} \sum_{m,l\ge1} D(N,m,l) x^N y^m z^l = \f{ U_0(x,y) U_0(x,z) }{ U_0(x,0) }. 
\end{equation}
\end{lemma}
We give the proof of Lemma \ref{lemma.d} in the appendix.

The number of disk triangulations with boundary length $m$ and a $02$-edge is determined by 
$[y^{m-3} z] W(x,y,z)$. 
To correct \eqref{Adef1} we should subtract a sum over all invalid layers, 
\begin{eqnarray}
\sum_{L: (E)} c_{k-2} x^n y^{l-2} 
&=& \sum_{k=2}^\infty c_{k-2} \sum_{l_1\ge 3} (l_1-2) x y^{l_1-4} [t^{l_1-3} z] W(x,t,z) \nn\\
&& \times \B( x U_0(x,0)\B)^{k-1},
\label{Adef2}
\end{eqnarray}
We need two intermediate calculations,
\[ \sum_{k=2}^\infty c_{k-2}\B( x U_0(x,0)\B)^{k-1} = u(x,0) F(u(x,0)), \]
\[ [t^{l_1-3} z] W(x,t,z) = \f{U_0'(x,0)}{U_0(x,0)} [t^{l_1-3}] U_0(x,t) = 
   \f{u'(x,0)}{u(x,0)} [t^{l_1-3}] U_0(x,t), \]
then we can continue with \eqref{Adef2}:
\begin{eqnarray}
\ldots 
&=&  u'(x,0) F(u(x,0)) \sum_{l_1\ge 3} (l_1-2) y^{l_1-4} [t^{l_1-3}] u(x,t) \nn\\
&=&  u'(x,0) F(u(x,0)) y^{-1} \f{\d}{\d y} ( y u(x,y) ) \nn\\
&=&  y^{-1} \B(u(x,y)+yu'(x,y)\B) u'(x,0) F(u(x,0)). \label{Adef3}
\end{eqnarray}
By subtracting \eqref{Adef3} from \eqref{Adef1} we get \eqref{Adef}. 
The proof is finished.
\eop

As an consequence of Lemma \ref{lemma.a} and Lemma \ref{lemma.lseq} we get
\begin{lemma}\label{lemma.ar}
Let $l$ and $k$ denote the length of lower and upper boundaries of a hull $B_R$,
and let $n$ denote the number of triangles in $B_R$. 
Take a function $F(t)$ that allows expansion \[ F(t)=\suml_{j=0}^\infty c_j t^j. \]
Then we can write down the sum over all $R$-hulls
\[ \sum_{B_R} c_{k-2} x^n y^{l-2} = (A^R(x)F)(y). \]
\end{lemma}

To show how Lemma \ref{lemma.ar} could be used, let us consider a $R$-hull $\barB_R$ and compute the 
expectation of the upper boundary $k$ when the lower boundary $l$ is fixed.
By Theorem \ref{thm.limdist},
\[ \E(k|l) = \sum_{B_R(l)} k \P(B_R) = \sum_{B_R(l)} k \f{b(k)}{b(l)} x_0^n, \]
where the sum is over all layers with lower boundary $l$. Applying Lemma \ref{lemma.ar}, we get
\[ \E(k|l) = \f1{b(l)} [y^{l-2}] (\bA_0^R f_1)(t), \]
where $\bA_0=\bA(x_0)$, 
\[ f_1(t) = \sum_{k=2}^\infty k b(k) = \sum_{k=2}^\infty k [y^{k-2}] B(y) t^k
   = \f1y \f{\d}{\d y} (y^2 B(y)).  \]
In a similar way an appropriate function can be constructed to compute $\E(k^2|l)$, $\E(k^3|l)$ 
and so on. 

An important fact is that the sum of probabilities over layers with a fixed lower boundary $l$ 
equals one for all $l$. This fact can be used to prove that the limiting measure $\tau$ on 
the space of triangulations $\calT$ is a probability measure, thus giving an independent 
proof of Theorem \ref{thm.as}.
\begin{lemma}\label{lemma.sumbr}
\begin{equation}\label{eq.sumbr}
 \sum_{B_R} \P\{B_R\RN S(m_0)\} = 1.
\end{equation}
\end{lemma}
\proof Let $f_0=B$. We have to show that
\[ \sum_{B_R} \P\{B_R\RN S(m_0)\} 
   = \sum_{B_R} \f{b(m_1(B_R))}{b(m_0)} x_0^{n(B_R)}
   = \f{[t^{l-2}] (\bA_0^R f_0)(t) }{[t^{l-2}f_0(t)]}
   = 1.
\]
This follows from the equality
\[ \bA_0^R f_0 = f_0, \]
which in it's turn follows from $\bA_0 f_0 = f_0$. 
To check the last equation one has to perform a trivial yet cumbersome expression transform. 
We omit it.
\eop

\subsection{Iterating linear operator}

To obtain asymptotic of moments as $R\to\infty$ we have to compute $\bA_0^R$.
First, we do a change of variable: 
\begin{equation}\label{syw}
w=\sqrt{1-\sqrt6 y}.
\end{equation}
When switching to such a new "coordinate system", a function $f(y)$ turns into a function $g(w)$
defined by the equation
\[ f(y) = g(w), \]
and the operator $\bA_0$ turns into the operator $\bB$, such that for all $f$, $y$
\[ (\bA_0 f)(y) = (\bB g)(w). \]
The operator $\bB$ acts on $g$ as follows
\[
(\bB(G))(w) = \f{ (2w+1)g\left(\f{w}{1+w}\right) }{ (1+w)^5 (1-w^2) }
           -\f{ \f18 (w+2) g\left(\f12\right) }{ (1+w)^2 (1-w^2) }.
\]
We can decompose it in two parts $\bB=\bB_1-\bB_2$,
\begin{eqnarray*}
\bB_1(g)(w) &=& \f{ (2w+1)g\B(\f{w}{1+w}\B) }{ (1+w)^5 (1-w^2) }
- \f3{32}\f{g\B(\f12\B)}{w^3(1-w^2)},\\
\bB_2(g)(w) &=& \f{ \f18 (w+2) g\left(\f12\right) }{ (1+w)^2 (1-w^2) } 
- \f3{32}\f{g\B(\f12\B)}{w^3(1-w^2)}.
\end{eqnarray*}
The reason for such decomposition is that both $\bB_1$ and $\bB_2$ are \defined{iterable},
i.e. we can write a simple formula for the $k$th iteration of each operator.
Let 
\[ %\begin{equation}\label{eq.g2G}
  g(w) = \f{G(w)}{w^3 (1-w^2)},
\] %\end{equation}
then
\[ %\begin{equation}\label{B1iter}
\bB_1^k(g)(w) = \f{ G\left(\f{w}{kw+1}\right) }{w^3(1-w^2)}
               - \f{ G\left(\f{1}{k+1}\right) }{w^3(1-w^2)},
\] %\end{equation}
\begin{equation}\label{eq.B1_1}
\bB_1^k(g)(1) = - \f{G'\left(\f1{k+1}\right)}{2(k+1)^2},
\end{equation}
\begin{equation}\label{eq.B2}
\bB_2(g) = \phi(z) g(\f12), \qquad 
\phi(z) = \f18 \f{(w+2)}{(1+w)^2 (1-w^2) } - \f3{32}\f{1}{w^3(1-w^2)}.
\end{equation}
To compute $k$th iteration of $\bB$, introduce the \defined{generating operators}
\[ \bM = \bM(z) = \sum_{k=0}^\infty z^k \bB^k, \qquad
   \bM_1 = \bM_1(z) = \sum_{k=0}^\infty z^k \bB_1^k, \]
and use an equality
%\EQ{Meq}{ \bM = \bM_1 + z \bM_1 \bB_2 \bM. }
\[ \bM = \bM_1 + z \bM_1 \bB_2 \bM. \]
Then
\EQ{Meqgw}{ \bM(g)(w) = \bM_1(g)(w) + z \bM_1(\phi)(w) \cdot \bM(g)(\f12), }
at $w=1/2$ this gives
\EQ{Meqgw12}{ \bM(g)(\f12) = \f{\bM_1(g)(\f12)}{1-z\bM_1(\phi)(\f12)}. }
Substitute \eqref{Meqgw12} to \eqref{Meqgw}, we'll get
\EQ{Meqgw2}{ \bM(g)(w) = \bM_1(g)(w) + \bM_1(g)(\f12) \theta(z,w), }
\[ \theta(z,w) = \f{ z \bM_1(\phi)(w) }{ 1-z\bM_1(\phi)(\f12) }. \]
There is however a simpler expression for $\theta(z,w)$. Under a change of variable $y\to w$
the function $f_0(y)=(1-\sqrt6y)^{-3/2}$ becomes $g_0(w)=1/w^3$. Substitute $g_0$ 
into \eqref{Meqgw2}, then since $\bB$ keeps $g_0$ we'll get
\[ \bM(g_0)(w) = \f{ g(w) }{ 1-z }, \] 
consequently
\[ \theta(z,w) = \f{ \bM_1(g_0)(w) - \f{ g_0(w) }{ 1-z } }{ \bM_1(g_0)(\f12) }. \]
We can summarize the computations above in a single lemma.
\begin{lemma}\label{lemma.M}
\[ %\begin{equation}\label{Meq_1}
\bM(g)(1) = \bM_1(g)(1) - \bM_1(g)(\f12) \f{ H_1(z)  - \f1{1-z} }{ H_2(z) },
\] %\end{equation}
where
\[
H_1(z) = \f1z \sum_{k=0}^\infty \f{z^k}{(k+1)^3}, \quad
H_2(z) = \f{32}{3}\sum_{k=0}^\infty \B( \f1{(k+1)^2} - \f1{(k+2)^2} \B) z^k.
\]
\end{lemma}
Now we pass to limiting distribution of $R$-hull's boundary as $R\to\infty$.
Using \ref{lemma.M}, let us compute the asymptotic of hull's upper boundary moments.
\begin{lemma}
\[ f_j(y) = \sum_{k=0}^\infty k^j c_k y^k
          = \f{(2j+1)!!}{2^j} \f1{w^{2j+3}} P_j(w^2), \]
where $P_j$ is a polynomial such that $P_j(0)=1$.
\end{lemma}
\proof
Proof by induction. First,
\[ f_0(y) = \Phi(y) = g_0(w). \]
Since 
\[ f_{j+1}(y) = y \f{d}{dy} f_j(y) \] 
and according to \eqref{syw}
\[ \f{\d w}{\d y} = -\f{\sqrt6}{2w}, \]
then
\begin{eqnarray*}
f_{j+1}(y) 
&=& \f{(2j+1)!!}{2^j} \B(\f{1-w^2}{\sqrt6} \B)
    \f{d}{d w}\B( \f1{w^{2j+3}} P_j(w^2) \B) \B(-\f{\sqrt6}{2w} \B) \\
&=& \f{(2j+1)!!}{2^j} \B( \f{2j+3}{2} \f{P(w^2)}{w^{2j+5}} - \f{P'(w^2)}{w^{2j+3}}\B)
    (1-w^2) \\
&=& \f{(2j+3)!!}{2^{j+1}} \f{1}{w^{2j+5}} P_{j+1}(w^2), 
\end{eqnarray*}
where
\[
P_{j+1}(t) = \B(P_j(t) - \f{2tP'(t)}{2j+3}\B)(1-t), \quad
P_{j+1}(0) = 1. 
\]
The lemma is proved.
\eop

\begin{lemma}\label{lemma.Gj}
Let \[ g_j(w) = \f1{w^{2j+3}}.  \]
Then
\EQ{bmz}{ \bM(g_j)(1) = \f{(2j)!}{(1-z)^{2j+1}} + O\B( \f{1}{(1-z)^{2j+2}} \B), }
\EQ{zrm}{ [z^R] \bM(g_j)(1) = R^{2j} + O(R^{2j-1}). }
\end{lemma}
\proof
Let 
\[ G_j(w) = \f{1-w^2}{w^{2j}}. \]
Then according to \eqref{eq.B1_1}
\begin{eqnarray*}
\bB_1^k(g_j)(1) 
   &=& - \f{G_j'\left(\f1{k+1}\right)}{2(k+1)^2}
    =  j(k+1)^{2j-1} - (j-1)(k+1)^{2j-3}, \\
\bB_1^k(g_j)(\f12) 
   &=& \f{32}{3} G\left(\f1{(k+2)}\right) - \f{32}{3}G\left(\f1{(k+1)}\right) \\
   &=& \f{32}{3} \B((k+2)^{2j}-(k+2)^{2j-2} - (k+1)^{2j}+(k+1)^{2j-2} \B).
\end{eqnarray*}
Let
\[ s_j(z) = \sum_{k=0}^\infty (k+1)^j z^k. \]
Then
\[ s_0(z) = \f1{1-z}, \qquad s_{j+1}(z) = \f{d}{dz} z s_j(z) = z s_j'(z) + s_j(z), \]
consequently as $z\to1$
\[ s_j(z) = z^j \f{d^j}{d z^j} s_0(z) + \ldots 
          = \f{j!}{(1-z)^{j+1}} + O\B( \f1{(1-z)^j} \B), \]
\begin{eqnarray*}
\bM_1(g_j)(1) &=& j s_{2j-1} - (j-1) s_{2j-3} \approx \f{(2j)!}{2(1-z)^{2j}}, \\
\bM_1(g_j)(\f12) &=& \f{32}{3}\left(\f1z-1\right) \B(s_{2j}(z) - s_{2j-2}(z)\B)
\approx \f{32}3 \f{(2j)!}{(1-z)^{2j}}.
\end{eqnarray*}
Using lemma \ref{lemma.M} we'll get \eqref{bmz}:
\[ \bM(g_j)(1) \approx \f{(2j)!}{2(1-z)^{2j}} - \f{32}3 \f{(2j)!}{(1-z)^{2j}}
   \f{ \zeta(3) - \f1{1-z} }{ \f{32}3 } \approx \f{(2j)!}{(1-z)^{2j+1}}.
\]
From the last expression \eqref{zrm} immediately follows.  The lemma is proved.
\eop

Now we pass to the proof of \ref{thm.profile}. Let $m_1(B_R)$ be the length of upper boundary
of a $R$-hull. As a consequence of \ref{lemma.Gj} 
\[ \E_R m_1(B_R)^j = \f{(2j+1)!!}{2^j} R^{2j} + O(R^{2j-1}),  \]
and the limit exists 
\[ \lim_{R\to\infty} \E \B( \f{m_1(B_R)}{R^2} \B)^j = \E \xi^j = \f{(2j+1)!!}{2^j}. \]
The moments generating function of $\xi$ is
\[ \sum_{j=0}^\infty \E\xi^j \f{(-s)^j}{j!} = \f1{(1+s)^{3/2}}. \]
Applying reversed Laplace transform we'll get the density
\[ p_{\xi}(x) = \f{2}{\sqrt\pi}e^{-t}t^{1/2}. \]
Theorem \ref{thm.profile} is proved.

\section{Branching process}\label{section.bp}

\subsection{Modified branching process}\label{section.mbp}
Now we can formally specify the {\em branching process with special behavior near zero}, 
that appears in Theorem \ref{thm.bp}.

Let $\zeta(t)$ be a branching process with offspring generating function
\begin{equation}\label{eq.phi}
\phi(t) = 1 - \f1{\dd \B(1+\f1{\sqrt{1-t}}\B)^2}.
\end{equation}
Let $X_R$ be the set of all trajectories of such a branching process on time interval $[0,R]$,
for all values of $\zeta(0)$. Each element of $X_R$ is then a forest of rooted trees, 
with height not exceeding $R$.

It's natural to assign to each trajectory $x\in X$ it's weight $\omega(x)$, equals to the probability
for $x$ to be a trajectory of $\zeta$. $\omega(x)$ can be represented as a product over
all vertices below $R$-th level in $x$,
\begin{equation}\label{eq.phiprod}
\omega(x) = \prod_{ v\in x, h(v)<R } [t^{desc(v)}] \phi(t), 
\end{equation}
where $desc(v)$ denotes the number of descendants of a vertex $v$.

Now consider a modified generating function $\barphi$,
\begin{equation}\label{eq.barphi}
\barphi(t) = \phi(t) - \f16 t \phi(t) + \f16.
\end{equation}
The modified branching process $\barzeta$ is defined as follows. To each trajectory $x\in X_R$ 
it assigns weight $\baromega(x)$ that differs from $\omega(x)$ in a single case:
if a whole generation of a branching process is inherited from a single vertex $v$, 
then in a term of a product \eqref{eq.phiprod} corresponding to $v$, 
the function $\phi$ should be replaced with $\barphi$.

For a system of weights $\baromega$ to define properly a probability distribution on a 
set of trajectories with a fixed starting state, we also have to modify the probability of 
degeneration. I.e. the transition $k\to 0$ should have probability $\barphi(0)\phi^{k-1}(0)$
instead of $\phi^k(0)$.

\subsection{Proof of Theorem 5} %5 stands for \ref{thm.bp}
Let us compare \eqref{eq.phiprod} to \eqref{eq.Zprod}. 
Both expressions assign to a trajectory of a branching process some weight, determined 
by the number of descendants at each vertex. 
The terms assigned to a vertex with $d$ descendants in \eqref{eq.phiprod} and \eqref{eq.Zprod}
are nearly the same, except that the function $\phi$ is {\em normalized} to satisfy the 
generating function conditions. Namely
\[ [t^d]\phi(t) = y_0^{d-1} [y^d] x_0 U_0(x_0,y), \]
which follows from the identity
\[ \phi(t) = y_0^{-1} x_0 U_0(x_0,y_0 t). \] 
(to check it one should use \eqref{eq.u0} for $U_0$ and note that according to \eqref{eq.defh} 
$h(x_0)=\f1{\sqrt6}=y_0$, the rest is trivial).

Imagine a trajectory $x$ and write $y_0$ on the top and $y_0^{-1}$ on the bottom of each edge.
A vertex $v$ then gets a term $y_0^{d-1}$, so for a trajectory $x\in X_R$ with starting (top) 
state $k$ and final (bottom) state $l$ we get
\begin{equation}\label{eq.pp}
\prod_{v\in x, h(x)<R } x_0 [y^d(v)] U_0(x_0,y)
    = y_0^{l-k} \prod_{v\in x, h(x)<R} [t^d(v)] \phi(t).
\end{equation}
It's easy to check that the modification $\barphi$ reflects the exception (E):
\[
[t^d] \barphi(t) = y_0^{d-1} [y^d] \B( x_0 U_0(x_0,y) - [z] x_0 W(x_0,y,z) \B),
\]
(see proof of Lemma \ref{lemma.a}).

Another important difference between the raw approach in section \ref{section.raw} and the 
branching process approach is the position of a {\em coordinate system origin}.
In a skeleton (being considered as a multi-rooted triangulation) the root on the lower boundary 
can be placed in an arbitrary position, while the root on the upper boundary is determined
via the root propagation rules (see layer pattern definition at page \pageref{def.lp}).
On the contrary, for the branching process trajectory (being drawn as a planar tree) we assume 
some order of particles at starting time; the order of particles for all subsequent generations
including the final one is induced.

Let $X_1$ be the set of all skeletons-trajectories with a root specified on the lower boundary only,
$X_2$ -- with the root on the upper boundary only, and $X_{12}$ -- with two roots specified
arbitrary on both boundaries. For an element $x\in X_1$ with upper boundary $k$ there are
$k$ elements in $X_{12}$; for and element $x\in X_2$ with lower boundary $l$ there are $l$
elements in $X_{12}$. 
Then for any function $f$ on skeletons that does not depend on the root position 
\begin{equation}\label{eq.fkl}
\sum_{x\in X_1} f(x) = \sum_{x\in X_{12}} k(x) f(x) = \sum_{x\in X_1} \f{k(x)}{l(x)} f(x).
\end{equation}
Thus, from \eqref{eq.Zprod}, \eqref{eq.pp} and \eqref{eq.fkl}
\[
\P\{ m_1(B_R)=k | m_0(B_R)=l \} 
 = \f{b(k)y_0^{-k}/k}{b(l)y_0^{-l}/l} \P\{\barzeta(0)=l|\barzeta(-R)=k\}. 
\]
To finish the proof of Theorem \ref{thm.bp}, note that 
\[ \sum_{k=0}^\infty  b(k)y_0^{-k}t^k/k = \int_{0}^t B(\f{y_0^{-1} \theta})\,d\theta  
 = \f2{\sqrt{1-t}} -4 + 2\sqrt{1-t}. \]

\subsection{Linear contour}
The main idea in establishing linear contour existence is counting ancestors.
The figure \figref{zz} shows, how a zigzag path in a skeleton that is allowed to bounce 
between two levels can be shorter that a horizontal one. 
\putfigure{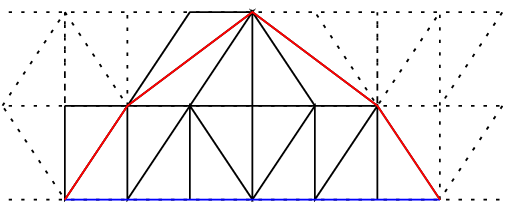}{A segment of a zigzag path (red)}

Given $r$th level with $n=xr^2$ edges, consider the corresponding $(-r)$th generation 
of a branching process, and count it's ancestors in the $(-2r)$th generation.
If this number is finite, there is a zigzag contour between levels $r$ and $2r$ with 
finite number of parts, i.e. it's length is linear in $r$. 

For such an estimate use the unmodified branching process $\zeta$. Since for each trajectory
$x$ it's unmodified weight majorises the modified weight, $\omega(x)\ge\baromega(x)$,
the expectation of any positive function with respect to $\omega$ estimates the expectation with respect to $\baromega$. 

Another reason to prefer $\zeta$ is that it's generating function is easily iterable, 
the $r$th iteration of $\phi$ equals
\begin{equation}\label{eq.iterphi}
 \phi_r(t) = 1+\f1{\dd\B(r+\f1{\sqrt{1-t}}\B)^2}. 
\end{equation}
For $\barzeta$ even computing the probability of non-degeneration in $r$th generation is 
a nontrivial task.

For a skeleton-trajectory $x$ denote by $a(x)$ the number of particles in moment $(-r)$
(upper boundary) that do have a non-empty offspring at moment $0$ (lower boundary)
(in other words this is the number of ancestors we wish to estimate).
Then an equality holds
\[
\sum_{x: k(x)=k} \omega(x) x^k t^l z^a = \B( z(\phi_r(t)-\phi_r(0)) + \phi_r(0) \B)^k.
\]
Consequently, the expected number of ancestors at level $2r$, conditioned to the level 
$r$ having $n$ vertices, is 
\[ \f{ [t^n]\XX(t) }{ [t^n] F_0(t)} , \]
where
\begin{eqnarray*}
\XX(t) &=&
  \f{\d}{\d z} F_0\B( z(\phi_r(t)-\phi_r(0)) + \phi_r(0) \B) \B|_{z=1} \\
&=& \vphantom{\B|} F_0'(\phi_r(t))(\phi_r(t)-\phi_r(0))
\end{eqnarray*}
and $F_0$ is modification of $f_0$ corresponding to the branching process approach,
\begin{eqnarray*}
F_0(t)  &=& \sum_{j=0}^\infty \f{b(j)}{j} y_0^{-j} t^j = \f2{\sqrt{1-t}} - 4 + 2\sqrt{1-t}, \\ 
F_0'(t) &=& \f{t}{(1-t)^{3/2}}.
\end{eqnarray*}
Using \eqref{eq.iterphi} for $\phi_r$ we get
\begin{eqnarray*}
\XX 
&=&  \f{ \phi_r(t) }{(1-\phi_r(t))^{3/2}}(\phi_r(t)-\phi_r(0))  \\
%&=& \B(1-\f1{\B(r+\f1{\sqrt{1-t}}\B)^2}\B) \B(\f1{r+1}^2-\f1{\B(r+\f1{\sqrt{1-t}}\B)^2}\B)
%    \B(r+\f1{\sqrt{1-t}}\B)^3 \\
&=& \f{1}{(r+1)^2} (1-t)^{-3/2}
    +\f{3r}{(r+1)^2} (1-t)^{-1}
    +\f{2(r^2-r-1)}{(r+1)^2} (1-t)^{-1/2}  \\
&&{}-\f{2r}{r+1}
    + \f1{r+(1-t)^{-1/2}}.  
\end{eqnarray*}
\begin{eqnarray*}
{} [t^n]\XX 
&=& \f{1}{(r+1)^2} \f{\G(n+1/2)}{\sqrt\pi n!}
    +\f{3r}{(r+1)^2} 
    +\f{2(r^2-r-1)}{(r+1)^2} \f{\G(n+3/2)}{\sqrt\pi n!}  \\
&&{}+[t^n]\f1{r+(1-t)^{-1/2}}. 
\end{eqnarray*}
\begin{eqnarray*}
\f1{r+(1-t)^{-1/2}} &=& \f{r-(1-t)^{-1/2}}{r^2-(1-t)^{-1}} \\
&=& \B( r-\sum_{j=0}^\infty \f{\G(j+\f12)}{\sqrt\pi j!} t^j \B)
    \B( \f1{r^2-1} + \sum_{j=1}^\infty \f{r^{2j-2}}{(r^2-1)^{j+1}} t^j \B), \\
{}[t^n] \f1{r+(1-t)^{-1/2}} &\le& \f{r^{2n-1}}{(r^2-1)^{n+1}}
        \approx \f{e^x}{r^3}.
\end{eqnarray*}
Finally,
\begin{eqnarray*}
[t^n] F_0(t) &=& \f{2(n-1)\G(n-\f12)}{\sqrt\pi n!} = \f2\pi n^{-1/2} + O(n^{-3/2}), \\
\f{ [t^n]\XX }{ [t^n]F_0(t) } &=& \f1{(r+1)^2} \f{4n^2-1}{4(n-1)}
   + \f{3r}{(r+1)^2} \B( \f{2(n-1)\G(n-\f12)}{\sqrt\pi n!} \B)^{-1} \\
&& + \f{2(r^2-r-1)}{(r+1)^2}  \f{2n-1}{4(n-1)} 
   + \f{ [t^n]\f1{r+(1-t)^{-1/2}} }{ [t^n] F_0(t) } \\
&\approx& x + \f{3\sqrt\pi}2 \sqrt{x} + 1 + O(\f1r).
\end{eqnarray*}
According to Theorem \ref{thm.profile}, the number of edges at $r$th level of a skeleton 
is approximately $\xi r^2$, so the expected length of a linear contour doesn't exceed 
\[
r \E \B(\xi+\f{3\sqrt\pi}2 \xi^{1/2} + 1 \B) < 10 r
\]
for large $r$.

\section{Universality}
Here we briefly discuss the universality of skeleton construction and the branching process approach.

In \cite{AS,A} two types of triangulations are considered.
Type II triangulations are the rooted near-triangulations, they are required to have no loops;
type III triangulations are \defined{strict near-triangulations}, they additionally required 
to have no double edges.
\begin{definition}
A \defined{strict near-triangulation} is a planar map with all faces being triangles
except the rooted face and with no cut vertices or double edges.
\end{definition}
For both types of triangulations the estimates for the $R$-hull volume and boundary length are
given and these values have the same order -- $R^2$ and $R^4$, up to polylogarithmic terms.

In present work we considered type II triangulations only.
However the skeleton construction is likely to work for type III triangulations too and lead
to a similar modified branching process.
\begin{conjecture}
For type III triangulations one can construct a branching process $\zeta_3$ and it's modification
$\barzeta_3$, so that the analog of Theorem \ref{thm.bp} holds. 
The branching process $\zeta_3$ is critical, has infinite variance and has non-degeneration probability
\[ \P\{ \zeta_3(r) > 0 | \zeta(0)=1 \} = \f1{r^2} + O(r^{-3}).  \]
\end{conjecture}
{\bf Motivation.}
The analogs of theorems \ref{thm.limdist} and \ref{thm.oneend} for type III triangulations exist. 
The ball and hull definition are the same. The skeleton can be defined the same way, 
since when gluing the slots of a skeleton with 
valid type III triangulations no double edges appears (in the general case).

The only things to change are the triangulations generating function (consequently, the 
offspring g.f. of a branching process) and the exception (E).
The exception should be formulated as follows:
$$
\begin{array}{ll}
\mbox{1)}&\mbox{all of $l_1,\ldots,l_k$ are equal to $2$, except $l_{j_1}$ and $l_{j_2}$; }  \\
\mbox{2)}&\mbox{both disks $D_{j_1}, D_{j_2}$ contain an edge between the vertices }\\
         &\mbox{$0$ and $2$ of the rooted boundary.}
\end{array}
\eqno (E3)
$$
Concerning the offspring generating function 
a preliminary computation show that it is the same for both $\zeta$ and $\zeta_3$,
i.e. the unmodified branching processes is the same for both types of triangulations.
However the {\em boundary coefficients}$b(m)$ 
for strict triangulations will be slightly different.

Another conjecture concerns the equivalence of modified and unmodified branching processes for 
large $r$ asymptotic.
\begin{conjecture}
For some class of functions $f$ (a.e. such that were used in section \ref{section.raw} for 
estimation of upper boundary moments asymptotic), the result of asymptotic estimation
with respect to modified or unmodified weight system/branching process is identical up 
to constant term, that doesn't depend on $f$:
\[
\lim_{R\to\infty} \f{ \sum_k c_k \P\{\zeta(0)=l|\zeta(-R)=k\} }
                      { \sum_k c_k \P\{\barzeta(0)=l|\zeta(-R)=k\} } = const(l).
\]
\end{conjecture}
{\bf Motivation.} 
Each trajectory $x$ can be broken into two parts by the first point, where the exception (E) is 
applicable, i.e. where the whole generation is inherited from a single parent.
Then for the upper part of trajectory $x_1$, $\omega(x_1)=\baromega(x_1)$,
while the lower part doesn't depend on $F$ and is likely to have finite length as $r\to\infty$.

The two in the limit sums above are then both parcels of sequences $a$, $b$ for $\zeta$ and 
$a$, $\bar b$ for $\barzeta$, 
\[
\f{ \suml_{r=0}^R a(R-r) b(r) }{ \suml_{r=0}^R a(R-r) {\bar b}(r) }.
\]
If both $b(r)$ and $\bar b(r)$ are decreasing as $1/r^2$ (which is likely to be true in our case)
and $a(r)$ grows as $r^\alpha$, $\alpha>2$, 
the limit of the expression above exists and depends on $b$ and $\bar b$ only.

\section*{Appendix.}

\proofof{Theorem \ref{thm.limdist}}
Since $T$ is rigid, 
\begin{eqnarray*}
\P\{ B\RN S_N \} &=& \f{\sum_{N_1+\ldots+N_k=N-n} C(N_1,m_1)\cdots C(N_k,m_k)}{C(N,m_0)} \\
&=& \f{ [x^{N-n}] \B( u_{m_1}(x) \cdots u_{m_k}(x) \B) }{ [x^N] u_{m_0}(x) },
\end{eqnarray*}
where
\[ u_m(x) = [y^{m-2}] U_0(x,y). \]
The behavior of the product $u_{m_1}(x)\cdots u_{m_k}(x)$ near singularities
is governed by the term $t^{3/2}$ of series 
\begin{eqnarray*}
\lefteqn{ u_{m_1} (x) \cdots u_{m_k}(x)\B|_{x=\pm(x_0-t)} = } &&\\ 
&& (\pm 1)^{m_1+\cdots+m_k} \prod_{j=1}^k a(m_j)  
\B[ 1 + \sum_{j=1}^k \f{b(m_j)}{a(m_j)} t 
+ \f{6^{3/4}}3 \sum_{j=1}^k \f{c(m_j)}{a(m_j)} t^{3/2} \B] + O(t^2) .
\end{eqnarray*}
We have to take into account both singularities $x=\pm x_0$. Thus 
\begin{eqnarray}
[x^{N-n}] u_{m_1}(x)\cdots u_{m_k}(x) &=& 
\f23 6^{3/4} \prod_{j=1}^k a(m_j) \cdot \sum_{j=1}^k \f{c(m_j)}{a(m_j)} \nn\\
&&{} \times [x^{N-n}](x_0-x)^{3/2} \B( 1+O(N^{-1/2} \B),
\label{uuxn}
\end{eqnarray}
as soon as $N-n = m_1+\cdots+m_k \mod 2$. This condition is satisfied, since 
$n=m_0+m_1+\cdots+m_k \mod 2$ and $N=m_0 \mod 2$.
In particular case of two boundaries, \eqref{uuxn} becomes
\begin{equation}\label{uuxn0}
[x^N] u_{m_1}(x) = \f23 6^{3/4} c(m_1) [x^N](x-x_0)^{3/2}\B(1+O(N^{-1/2})\B). 
\end{equation}
The statement of the theorem is then a fraction of \eqref{uuxn} and \eqref{uuxn0}.
Theorem \ref{thm.limdist} is proved.
\eop

\proofof{Theorem \ref{thm.oneend}}
Without loss of generality let $j=1$.
\[ \P\{ N_1=\max_{i=1,\ldots,k}(N_i) \} = 
   \f{ \suml_{N_1, \ldots, N_k \atop N_1\ge N_2,\ldots, N_k} C(N_1,m_1)\cdots C(N_k,m_k) }
     { \suml_{N_1, \ldots, N_k} C(N_1,m_1)\cdots C(N_k,m_k) }
\]
For each integer $M>0$ 
\begin{eqnarray}
                                                               \lefteqn{
\lim_{N\to\infty}                                             
                  \P\{ N_1=\max_{i=1,\ldots,k}(N_i) \} \ge 
\lim_{N\to\infty} \P\{ N_2+\ldots+N_k \le M \}                          } &&
\nn\\ 
&=& \suml_{N_2+\ldots+N_k \le M} C(N_2,m_2)\cdots C(N_k,m_k)
    \lim_{N\to\infty} 
    \f{ [x^{N-n-(N_2+\ldots+N_k)}] u_{m_1}(x) }{ [x^{N-n}] u_{m_1}(x) \cdots u_{m_k}(x) }
\nn\\
&=& \suml_{N_2+\ldots+N_k \le M} C(N_2,m_2)\cdots C(N_k,m_k) 
    \f{ c(m_1) x_0^{N_2+\ldots+N_k} }{\prodl_{i=1}^k a(m_i) \cdot \suml_{i=1}^k \f{c(m_i)}{a(m_i)}}.
\label{mp1}
\end{eqnarray}
Take a limit $M\to\infty$. We'll get
\begin{eqnarray}
\lim_{N\to\infty}                                         \lefteqn{
                  \P\{ N_1=\max_{i=1,\ldots,k}(N_i) \}            } &&
\nn\\
&\ge& \f{ c(m_1) }{ \prodl_{i=1}^k a(m_i) \cdot \suml_{i=1}^k \f{c(m_i)}{a(m_i)} }
      \lim_{M\to\infty} \suml_{N_2+\ldots+N_k \le M} \prodl_{i=2}^k C(N_i,m_i) x_0^i 
\nn\\
&\ge& \f{ c(m_1) }{ \prodl_{i=1}^k a(m_i) \cdot \suml_{i=1}^k \f{c(m_i)}{a(m_i)} }
      \lim_{M\to\infty} \prodl_{i=2}^k u_{m_i}(x_0) 
\nn\\
&=& \f{c(m_1)\prodl_{i=2}^k a(m_i) }{\prodl_{i=1}^k a(m_i) \cdot \suml_{i=1}^k \f{c(m_i)}{a(m_i)}}
 =   \f{c(m_1)}{a(m_1)} \B( \suml_{i=1}^k \f{c(m_i)}{a(m_i)} \B)^{-1}.
\label{mp2}
\end{eqnarray}
Thus for each $j=1,\ldots,k$
\begin{equation}
\label{ncam0}
\lim_{N\to\infty} \P\{ N_j=\max_{i=1,\ldots,k}(N_i) \} \ge 
 \f{c(m_j)}{a(m_j)} \B( \suml_{i=1}^k \f{c(m_i)}{a(m_i)} \B)^{-1}.
\end{equation}
Since the sum of right hand sides of \eqref{ncam0} for $j=1,\ldots,k$ is $1$, 
the inequality can be replaced by a strict equality. 
This gives the first statement of the theorem.

From \eqref{mp1}, \eqref{mp2} 
\[ \lim_{N\to\infty} \P\{N_1=\max_{i=1,\ldots,k}(N_i)\} = \lim_{M\to\infty}
   \lim_{N\to\infty} \P\{N_2+\ldots+N_k \le M \},
\]
and the limiting conditional distribution of 
$(N_2,\ldots,N_k|N_1=\max)$ has a generating function  
\begin{equation}\label{en1}
\lim_{N\to\infty} \E(t_2^{N_2}\cdots t_k^{N_k}|N_1=\max_{i=1,\cdots,k}N_i)
   =\prodl_{i=2}^k \f{u_{m_k}(x_0 t_i)}{ u_{m_k}(x_0) }.
\end{equation}
%Thus the above conditional distribution exists, 
Thus the conditional distribution $(N_2,\ldots,N_k|N_1=\max_{i=1,\ldots,k} N_i)$ exists
and the random variables $N_2,\ldots,N_k$ are asymptotically independent. 
Theorem \ref{thm.oneend} is proved.
\eop

\proofof{Lemma \ref{lemma.d}}
Let $T$ be a triangulation counted by $D(N,m,l)$, i.e. with a cut edge between vertices $0$ and $l$
of the boundary ($0l$-edge). 
Cut $T$ in two parts along this edge, if there are multiple $0l$-edges, choose the rightmost one
(i.e. the one that is met first when walking around the root vertex counter-clockwise starting 
from the root, see \figref{cut}).
Then one part $T_1$ has type $(N_1,l+1,0)$ and has no edge parallel to the root,
the second part $T_2$ has type $(N_2,m+1,0)$ and is a generic triangulation.

\putfigure{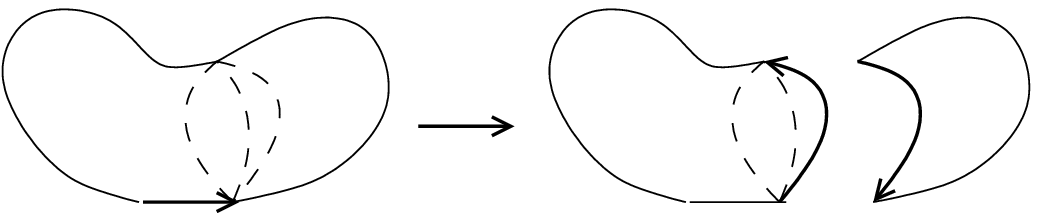}{Cutting the disk triangulation}

Let $R(N,m)$ be the number of triangulations of type $(N,m,0)$ with no edge parallel to the root,
and let $R(x,y)$ be the corresponding generating function
\[ R(x,y) = \sum_{N\ge0} \sum_{m\ge2} R(N,m) x^N y^{m-2}.  \]
Then 
\begin{equation}\label{eq.wur1}
W(x,y,z) = U_0(x,y) R(x,z).
\end{equation}
By definition $W(x,y,z)$ is symmetric in $y$, $z$ and $W(x,y,0)=U_0(x,y)$, consequently
\begin{equation}\label{eq.wur2}
U_0(x,y) = U_0(x,0) R(x,y).
\end{equation}
From \eqref{eq.wur1}, \eqref{eq.wur2} the statement immediately follows,
\[
W(x,y,z) = \f{ U_0(x,y) U_0(x,z) }{ U_0(x,0) }.
\]
\eop

%\vskip 3cm
%\noindent
%\begingroup
%\small
%Maxim Krikun\\
%Laboratory of Large Random Systems,\\
%Faculty of Mechanics and Mathematics,\\
%Moscow State University.\\
%Russia, Moscow, 119922, Vorobjevy Gory, MSU.\\
%{krikun@lbss.math.msu.su}\\
%\endgroup
%

\end{document}